\documentclass[12pt]{amsart}
\newcounter{defcounter}
\usepackage{amssymb,amsmath,amscd,graphicx, enumitem,
latexsym,amsthm}
\usepackage{amssymb,latexsym,eufrak,amsmath,amscd,graphicx, amsxtra}
  \usepackage[all]{xy}
  \usepackage{pdfsync}
  \usepackage[mathscr]{euscript}

\theoremstyle{plain}
\newtheorem{theorem}{Theorem}

\newtheorem{proposition}[theorem]{Proposition}

\newtheorem{lemma}[theorem]{Lemma}

\newtheorem*{propositionn}{Proposition}

\newtheorem{proposition.definition}[theorem]{Proposition/Definition}

\newtheorem{theoremalpha}{Theorem}

\theoremstyle{definition}
\newtheorem{definition}[theorem]{Definition}

\newtheorem{remark}[theorem]{Remark}
\newtheorem{example}[theorem]{Example}

\newcommand{\noi}{\noindent}

\newcommand{\rk} {\text{rank }}

\newcommand{\mult}{\textnormal{mult}}

\newcommand{\gon}{\textnormal{gon}}

\newcommand{\Hilb}{\textnormal{Hilb}}

\newcommand{\Sym}{\textnormal{Sym}}

\newcommand{\BVA}[1]{\textnormal{(BVA)}_{#1}}
\newcommand{\stabcorr}{\textnormal{stab.corr}}
\newcommand{\stabirr}{\textnormal{stab.irr}}
\newcommand{\irr}{\textnormal{irr}} 
\newcommand{\birat}{\textnormal{birat}}
\newcommand{\pr}{\textnormal{pr}}
\newcommand{\corr}{\textnormal{corr}}
\newcommand{\id}{\textnormal{id}}
\newcommand{\pt}{\textnormal{pt}}
\newcommand{\Prj}{\mathbb P}
\newcommand{\im}{\text{im}}

\numberwithin{theorem}{section}

\begin{document}

\title[Stabilized measures of association]{Stabilized measures of association: comparing the birational geometry of varieties having different dimensions}
 \author{Giovanni Passeri}
 \address{Department of Mathematics, State University of New York at Stony Brook, 100 Nicolls Rd, Stony Brook, NY 11794, USA}
\email{giovanni.passeri@stonybrook.edu}

\maketitle

\subsubsection*{Abstract} We extend some of the measures of association defined by Lazarsfeld and Martin, obtaining useful invariants to compare the birational geometry of two varieties having different dimensions. We explore such invariants providing examples and computing them in some cases of interest.

 \section*{Introduction}
\setcounter{tocdepth}{1}
In recent years there has been interest in studying \textit{measures of irrationality}, namely numerical measures of the failure for a variety to be rational: see for example \cite{BDELU}. 
More recently Lazarsfeld and Martin \cite{lazarsfeld2021measures} defined what they called \textit{measures of association}, which measure the failure for two varieties $X,Y$ of the same dimension to be birationally equivalent. The purpose of this note is to explore such measures in a context in which we allow the varieties to have different dimensions, defining what we will name, for lack of better terminology, \textit{stabilized measures of association}. As the name suggests, the idea is simply to consider products with projective spaces, to make the dimensions coincide. 
\vspace{0.25cm}

Turning to details, given two varieties $V_1,V_2$ of dimension $N$, for us a \textit{correspondence} between them is a an $N$-dimensional subvariety $Z \subset V_1 \times V_2$ dominating both $V_1$ and $V_2$; Lazarsfeld and Martin \cite{lazarsfeld2021measures} define the \textit{correspondence degree} between $V_1$ and $V_2$ to be 
\[
\text{corr}(V_1,V_2):=\min \Big\{\deg (Z \to V_1) \cdot \deg (Z \to V_2) \ \Big| \ Z\text{ is a correspondence between $V_1$ and $V_2$}\Big\}.
\] 
Given now two varieties $X,Y$ of dimensions $m \le n$ respectively, we will focus on two integers:
\[
\begin{aligned}
    & \text{corr}(X \times \mathbb P^{n-m}, Y) \text{ and }\\
    &\text{stab.corr}(X,Y):=\min_{m+r=n+s}\text{corr}(X \times \mathbb P^r, Y \times \mathbb P^s).
\end{aligned}
\]
We call the latter \textit{stable correspondence degree} between $X$ and $Y$. Note that 
\[
\stabcorr (X,Y)=1 \iff {
\begin{aligned}
    &\text{ $X$ and $Y$ are \textit{stably birational equivalent}}\\
    &\text{i.e.: $X \times \mathbb P^r \sim_{\birat} Y \times \mathbb P^s$ for some $r,s \ge 0$.} 
    \end{aligned}
}
\]

Taking fibered products we get

\begin{equation}\label{obvioous bound}
\text{stab.corr}(X,Y) \le \text{stab.irr}(X)\cdot \text{stab.irr}(Y).\footnote{Recall that the \textit{stable degree of irrationality} is an invariant defined by Bastianelli \cite{bastianelli2017irrationality} as $\text{stab.irr}(X):=\min_{r \ge 0}\text{irr}(X \times \mathbb P^r)$, where the \textit{degree of irrationality} $\text{irr}(V)$ of some variety $V$ is defined to be the least degree of a dominant rational map $V \dashrightarrow \mathbb P^{\dim V}$.}
\end{equation}
Moreover
\[
\log (\stabcorr)
\]
defines a metric on the set of projective varieties modulo stable birational equivalence.  

We think of the failure of equality in (\ref{obvioous bound}) as suggesting there are ``unexpected" correspondences $Z\subset (X \times \mathbb P^r) \times (Y \times \mathbb P^s)$. In particular we might expect equality in (\ref{obvioous bound}) if $X$ and $Y$ are chosen to be general in some suitable sense. Our main effort here is to make this intuition precise in several cases of interest.

To begin with we have the following proposition.

\begin{propositionn} Let $C$ and $D$ be very general curves of genera $g_C$ and $g_D$ respectively. Then
\[
\stabcorr(C,D)=\gon (C) \cdot \gon(D) = \stabirr(C) \cdot \stabirr(D).
\]
\end{propositionn}
\noi See Example \ref{ex curve vs curve} for a proof.

Our first main result asserts that equality holds in (\ref{obvioous bound}) for very general hypersurfaces of sufficiently large degree.

\begin{theoremalpha}\label{curves & hypersf}
Let $X\subset \mathbb P^{m+1}, Y \subset \mathbb P^{n+1}$ 
be very general hypersurfaces with degrees respectively $\ge 2m+2$ and $\ge 2n+2$. Then
\[
\textnormal{stab.corr}(X,Y) = \textnormal{stab.irr}(X) \cdot \textnormal{stab.irr}(Y).
\]
\end{theoremalpha}
\noi The statement of Theorem \ref{curves & hypersf} still holds if we replace $Y$ with a very general curve, see Theorem \ref{curves vs hypersf 2}. 

Our second result compares curves with abelian surfaces.

\begin{theoremalpha}\label{curves vs abelian surface} Let $C$ be a very general curve of genus $g>1$ and $A$ a very general abelian surface. Then
\[
\textnormal{corr} (C \times \mathbb P^1, A)= \textnormal{gon} (C) \cdot \irr (A).
\]

\end{theoremalpha}
 \noi Note that the degree of irrationality of a very general $(1,d)$-polarized abelian surface is not known for $d=1,3$, while the theorem holds also in those cases. We also remark that, even though Theorem \ref{curves vs abelian surface} holds true, in most cases there are correspondences $Z\subset (C \times \mathbb P^1) \times A$ with 
 \[
 \deg (Z \to A)<\irr (A).
 \]
 This is still true even if we replace $C \times \mathbb P^1$ with any surface, see Proposition \ref{p strong principle of generality fails for C,A} and Example \ref{ strong princ og gen fails for C, A}. The case $g=1$ remains open as does the question whether $\stabcorr (C,A)=\stabirr (C) \cdot \stabirr(A)$. 
 In fact, to the best of our knowledge, it is not even known yet what $\stabirr(A)$ is. However, we show that $\stabcorr(C,A)\ge 3 \cdot \gon (C)$, under the hypothesis of Theorem \ref{curves vs abelian surface}, see Proposition \ref{stab.corr(C,A)}.
Concerning organization: in $\S$\ref{section 1} we make some general remarks on the invariant we propose. In $\S$\ref{section 2} we present some examples of interesting correspondences. $\S$\ref{section 3} is devoted to the study of correspondences involving hypersurfaces of high degree, proving in particular Theorem \ref{curves & hypersf}. Similarly in $\S$\ref{section 4} we study correspondences involving curves and abelian surfaces, proving Theorem \ref{curves vs abelian surface}. 

We will work over the complex numbers throughout all the paper.

\subsection*{Acknowledgments} I would like to thank my advisor Robert Lazarsfeld for suggesting the problem and for his encouragement and guidance throughout the work on this project. I would also like to thank Olivier Martin for his helpful comments and for kindly reviewing a preliminary version of this note. Furthermore I would like to thank Nathan Chen, Prabhat Devkota, Lisa Marquand, Saverio Andrea Secci, Jason Starr and particularly Srijan Ghosh and Federico Moretti for engaging in valuable discussions.

\section{First remarks and examples}\label{section 1}
\numberwithin{equation}{section}
\setcounter{equation}{0}

We devote this section to some preliminary remarks and examples. We start by fixing notation. Let $V_1,V_2$ be two projective varieties of the same dimension. A \textit{correspondence} between them is given by a projective variety $Z$ having two generically finite dominant maps 
\[
a_1:Z \to V_1, \  a_2:Z \to V_2.
\]
and we define its \textit{degree} to be
\[
\deg Z:= \deg (a_1) \cdot \deg (a_2).
\]
Let now $Z$ be such a correspondence. Then $Z':= (a_1\times a_2)(Z)\subset V_1\times V_2$ is a correspondence and 
\[
\begin{aligned}
    &\deg (Z' \to V_1) \le \deg a_1\\
    &\deg (Z' \to V_2) \le \deg a_2.
\end{aligned}
\]
So, for the purpose of studying the minimal degree of a correspondence between $V_1$ and $V_2$, it is not restrictive to limit ourselves to correspondences $Z\subset V_1 \times V_2$ with maps given by the projections. So for the rest of the paper we will stick with this second notion. Then, we can view the fibers of $a_2$ (resp $a_1$) as cycles in $V_1$ (resp $V_2$). 
With this being said we fix the following set-up for the rest of the paper.

\subsection*{Set-Up}\label{general setup} Let $X,Y$ be any two irreducible projective varieties, of dimension $m$ and $n$ respectively. We consider a correspondence 
\[
Z\subset (X \times \mathbb P^r) \times (Y \times \mathbb P^s)
\]
between $X \times \mathbb P^r$ and $Y\times \mathbb P^s$, where $r,s\ge 0$ are integers such that $m+r=n+s$. We denote the projection maps 
\[
a:= \pr_{X \times \mathbb P^r}|_Z: Z \to X \times \mathbb P^r\text{ and }b:=\pr_{Y \times \mathbb P^s}|_Z: Z \to Y \times \mathbb P^s.
\]
We summarize the situation in the following diagram
\begin{equation}
\xymatrix{ & Z \ar[dr]^b \ar[dl]_a& \\
X \times \mathbb P^r & & Y \times \mathbb P^s}
\end{equation}

\subsection{Correspondences between varieties with canonical forms} 
In earlier works on measures of irrationality \cite{bastianelli2011gonality} \cite{BDELU} \cite{lopez1995curves}, a natural strategy has been to exploit positivity properties of the canonical bundle together with some vanishing statements. What we will do is to review and adapt the techniques used in \cite{lazarsfeld2021measures} to the case in which we allow varieties of different dimensions.
\vspace{0.5cm}
\\
Let $X,Y,Z$ be as in the Set-Up. Then, for every $k>0$, $Z$ defines homomorphisms
\[
\begin{aligned}
&Z_*^{k,0}:=\text{Tr}_b \circ a^* \circ \pr_X^*:H^{k,0}(X)\to H^{k,0}(Y \times \mathbb P^s)\\
&{Z^*}^{k,0}:=\text{Tr}_a \circ b^* \circ \pr_Y^*:H^{k,0}(Y)\to H^{k,0}(X\times \mathbb P^r),
\end{aligned}
\]
where $\pr_X: X \times \mathbb P^r \to X$ and $\pr_Y:Y \times \mathbb P^s \to Y$ are the projections and $\text{Tr}$ is the trace of holomorphic $k$-forms. We can write this more explicitly: if $\omega$ is a holomorphic $k$-form on $X$ and $(y,v)\in Y \times \mathbb P^s$ is a general point, we have
\begin{equation}\label{trace}
Z_*^{k,0}\omega(y,v)=\sum_{(x,u)\in b^{-1}(y,v)}\pr_X^*\omega (x,u)=\sum_{(x,u)\in b^{-1}(y,v)}\omega(x). \footnote{Note that in (\ref{trace}) we are identifying $T_{(x,u)}Z$ with $T_{(y,v)}(Y \times \mathbb P^s)$ via $db_{(x,u)}$, for all $(x,u)\in b^{-1}(y,v)$; this can be done for $(y,v) \in Y \times \mathbb P^s$ general, since $b$ is generically finite. Moreover with abuse of notation we are writing $\omega(x)$ instead of $\omega(x)({d\text{pr}_X}_{(x,u)}-)$.}
\end{equation}
So if we can find for some $(y,v) \in Y \times \mathbb P^s$ a $(k,0)$-form which fails to vanish at exactly one point of $\pr_X(b^{-1}(y,v))$, then $Z_*^{k,0}\ne 0$. The same considerations hold for ${Z^*}^{k,0}$. This brings into the picture the birational positivity of $K_X$ (resp $K_Y$).

Let $X$ be a projective algebraic variety. We say a line bundle $L$ on $X$ \textit{birationally separates $p+1$ points} or \textit{satisfies the property} $\BVA p$ if there is an open Zariski set $U \subset X$ such that for every $0$-dimensional subscheme $\xi \subset U$ of length $p+1$ the restriction map 
\[
H^0(L) \to H^0(L \otimes \mathcal O_\xi) \text{ is surjective}.
\]
In particular, we ask that for any choice of $p+1$ distincted points $x_0, \dots, x_p\in U$ we can find a section $s\in H^0(L)$ vanishing at $x_1, \dots x_p$ and not vanishing at $x_0$. If this is the case we say that $L$ \textit{separates} $\xi=x_0+ \dots +x_p$.
This leads to
\begin{proposition}
        Let $X,Y,Z$ as in the Set-Up. Suppose $K_X$ (resp $K_Y$) birationally separates $p+1$ points. If, for a general $(y,v)\in Y \times \mathbb P^s$ we have $\deg \pr_X (b^{-1}(y,v)) \le p+1$ (resp $\deg \pr_Y(a^{-1}(x,u)) \le p+1$ for a general $(x,u)\in X \times \mathbb P^r$), then the map
    \[
    Z_*^{m,0}:H^{m ,0}(X) \to H^{m,0}(Y \times \mathbb P^s)
    \]
    is non-zero (resp the map ${Z^*}^{n,0}:H^{n,0}(Y) \to H^{n,0}(X \times \mathbb P^r)$ is non-zero).
\end{proposition}

\begin{proof}
    Let $(y,v)\in Y\times \mathbb P^s$ be a general point. By hypothesis we can find an $(m,0)$-form $\omega$ on $X$ vanishing at all but one point of $\pr_{X}(b^{-1}(y,v))$. Call such a point $x_0$. Then
    \[
    \begin{aligned}
    &(Z_*^{m,0}\omega)(y,v)=\sum_{(x,u)\in b^{-1}(y,v)}\omega(x)=\omega(x_0)+\sum_{x\in \pr_{X}(b^{-1}(y,v)), \ x \ne x_0}\omega(x)=\omega (x_0)\ne 0,\\
    &\text{so $(Z_*^{m,0})\omega \ne 0$.}
    \end{aligned}
    \]
\end{proof}

\begin{example}\label{ex hypersf vs hypersf}
    Let $X\subset \mathbb P^{m+1}$, $Y\subset \mathbb P^{n+1}$ be very general hypersurfaces of degrees respectively $d \ge 2m+2, e \ge 2n+2$. Then, as in \cite{lazarsfeld2021measures} equation (1.7), one has
    \begin{equation}\label{morph HS hypersf}
    \begin{aligned}
        &\text{Hom}_\text{HS}(H_\text{prim}^m(X, \mathbb Z),H^m(Y, \mathbb Z))=0\\
        &\text{Hom}_\text{HS}(H_\text{prim}^n(Y, \mathbb Z),H^n(X, \mathbb Z))=0
    \end{aligned}
    \end{equation}
    See also \cite{PeterSteenbrinkMixedHS} Corollary 20.23.
    Let now $Z$ be as in the Set-Up. Consider the natural morphism of Hodge structures 
    \[
    p_Y:H^m(Y \times \mathbb P^s, \mathbb Z) \to H^m(Y, \mathbb Z).
    \]
    Since $H^{m,0}(X)\subset H_\text{prim}^m(X, \mathbb C)$, we have $p_Y\circ Z_*^{m,0}=0$ by example \ref{morph HS hypersf}. Moreover, $p_Y$ is an isomorphism at the $(m,0)$-level, so $Z_{m,0}^*=0$. Thus, for a general $(y,v)\in Y \times \mathbb P^s$ and for all $\omega \in H^0(K_X)$, we have
    \[
    0=(Z_*^{m,0}\omega)(y,v)=\sum_{x\in \pr_{X}(b^{-1}(y,v))}\omega(x).
    \]
Therefore $K_X$ does not separate $\pr_X(b^{-1}(y,v))$. Since $K_X=\mathcal O_X(d-m-2)$, by Lemma \ref{lemma Caley-Bacharach} (below) we have that that 
\begin{itemize}
    \item $\deg \pr_X (b^{-1}(y,v)) \ge d-m$ and
    \item if $\deg \pr_X (b^{-1}(y,v)) \le d-1$, then $\pr_X (b^{-1}(y,v))$ is made up of collinear points.
\end{itemize}
Similarly $ \deg \pr_Y(a^{-1}(x,u)) \ge e-n$ and if $\deg \pr_Y(a^{-1}(x,u)) \le e-1$, then $\pr_Y(a^{-1}(x,u))$ is made up by collinear points.
\end{example}

\begin{lemma}[\cite{bastianelli2011gonality} Lemma 2.4]\label{lemma Caley-Bacharach} Fix $N  \ge 2$ and let $W\subset \mathbb P^N$ be a set of $k \ge 1$ distincted points. Let $r\ge 1$ and suppose that no section of $\mathcal O_{\mathbb P^N}(r)$ separates $W$.
Then $k \ge r+2$. Moreover if $k \le 2r+1$, then $W$ lies on a line.
\end{lemma}

\begin{example}\label{ex curve vs curve}
    Let $X,Y$ be very general curves of genera $g_X,g_Y>0$. As in \cite{lazarsfeld2021measures} equation (1.6), we have
    \[
    \text{Hom}_\text{HS}(H^1(X, \mathbb Z), H^1(Y,\mathbb Z))=0
    \]
    see also   \cite{CilibertoJacobians} Theorem 3.1.
    Hence, if $Z$ is as in the Set-Up above, we have that $\pr_X(b^{-1}(y,v))$ is not separated by $K_X$, for general $(y,v)\in Y \times \mathbb P^s$. Therefore
    \begin{itemize}
        \item $\deg b \ge \deg \pr_X(b^{-1}(y,v)) \ge \gon (X) \text{ and similarly}$
        \item $\deg a \ge \gon (Y).$
    \end{itemize}
\end{example}
\begin{example}\label{ex curve vs hypersf}
    Let $X$ be any curve of genus $g_X \ge 0$ and $Y\subset \mathbb P^{n+1}$ be a very general hypersurface of degree $e\ge 2n+2$.
    As in Examples \ref{ex hypersf vs hypersf} and \ref{ex curve vs curve} we have: if $n=1$, then $\text{Hom}_\text{HS}(H^1(X, \mathbb Z), H^1(Y,\mathbb Z))=0$, while if $n>1$, $H^1(Y,\mathbb Z)=0$ by the Lefschetz Hyperplane Theorem. So in any case, both $Z_*^{1,0}$ and ${Z^*}^{n,0}$ vanish. So we get 
    \begin{itemize}
        \item $\pr_X(b^{-1}(y,v))$ is not separated by $K_X$, for a general $(y,v)\in Y \times \mathbb P^s$, so $\deg a \ge \gon (X)$;
        \item For $(x,u)\in X \times \mathbb P^r$ general, $\deg \pr_Y(a^{-1}(x,u)) \ge e-n$ and if $\deg \pr_Y(a^{-1}(x,u)) \le e-1$, then $\deg \pr_Y(a^{-1}(x,u))$ is made up of collinear points.
    \end{itemize}
\end{example}

\begin{example}\label{ex curve ab surf}
    Analogously, let $X$ be a very general curve and $Y$ a very general abelian variety. If $Z$ is as in the Set-Up we have
    $\deg a \ge \gon (X)$, with the same argument as in the previous examples.
\end{example}
\vspace{0.5cm}
\section{Examples of correspondences}\label{section 2}

In this section we will present examples of some interesting correspondences. We felt that presenting these early on could help bring to life the new phenomena that occur for stabilized correspondences.

First, let $X,Y$ as in Theorem \ref{curves & hypersf}. We will see in the proof of Theorem \ref{curves & hypersf} that, if $Z$ is as in the Set-Up, then 
\begin{equation}\label{principle of generality}
\deg a \ge \stabirr(Y) \text{ and }\deg b \ge \stabirr(X). 
\end{equation}
Examples \ref{bad curve and hypersf} and \ref{general bad hypersf} show this fails if $X,Y$ are not very general. 

Turning to Theorem \ref{curves vs abelian surface}, Example \ref{bad curve and ab surface} shows it fails if we  drop the ``very general" hypothesis. 

\noi
Moreover, even though Theorem \ref{curves vs abelian surface} ensures $\corr (C \times \mathbb P^1, A) = \gon (C) \cdot \irr (A)$ for a very general curve of genus $g>1$ and a very general abelian surface $A$, in most cases there are correspondences $Z\subset (C \times \mathbb P^1) \times A$ having 
\[
\deg (Z \to C \times \mathbb P^1)<\irr (A).
\]
In fact in Example \ref{ strong princ og gen fails for C, A} we will construct for any surface $S$ and for a very general abelian surface $A$, a correspondence $Z \subset S\times A$ with $\deg (Z \to S)=3$, while if $A$ is $(1,d)$-polarized with $d\nmid 6$ we have $\irr (A)=4$, see \cite{martin}.  

\begin{example}\label{bad curve and hypersf}
    Let $S\subset \mathbb P^3$ be a smooth surface of degree $d$ and $C\subset S$ be a smooth curve, e.g. a general hyperplane section of $S$. Consider the projective bundle $V:=\mathbb P_\text{sub}(T_S|_C)\to C$. A point in $V$ consists of a point $x\in C$ and a line $\ell_x\subset \mathbb P^3$ tangent to $S$ at $x$. Consider then the correspondence 
    \[
    Z:=\overline{\{(\ell_x,y) \in V \times S \ | \ y\in \ell_x, y \ne x\}}.
    \] 
    with projections $a:Z \to V, b: Z \to S$. Then for general $\ell_x\in V$,
    \[
    a^{-1}(\ell_x)=(\ell_x \cdot X)\setminus \{x \} \implies \deg a=d-2
    \]
    So in particural $\deg (Z \to V)<d-1=\stabirr(S)$, provided $d\ge 6$, see \cite{yang2019irrationality}. Since $V\sim_\text{birat}C \times \mathbb P^1$, this gives an example of a correspondence between $C \times \mathbb P^1$ and $S$ having degree over $C \times \mathbb P^1$ smaller than $\stabirr (S)$. 
\end{example}
\begin{example}\label{general bad hypersf}
    Let $X \subset \mathbb P^{m+1}$ be a hypersurface of degree $d$ and let $Y$ some variety of dimension $n$ with a rational map $f: Y \dashrightarrow X$ which is generically finite onto its image. Consider the map given by the composition
    \[
    \xymatrix{ Y \times \mathbb P^{m-n} \ar@{-->}[r]^i & Y \times \mathbb P^{m-1} \ar@{-->}[r]^{f\times \id}& X \times \mathbb P^{m-1} \ar@{-->}[r]^j & \mathbb P_\text{sub}(T_X)}
    \]
    where $i$ is some map birational onto its image and $j$ some birational equivalence over $X$. As in Example \ref{bad curve and hypersf}, the image of a general point point $(y,t)$ is some line $\ell_{(y,t)} \subset \mathbb P^{m+1}$ which is tangent to $X$ at $f(y,t)$. Then we can consider the correspondence
    \[
    Z:= \overline{\{(y,t,x) \in Y \times \mathbb P^{m-n} \times X \ | \ x\in \ell_{(y,t)}, x \ne f(y,t)\}}.
    \]
    As before this is a correspondence between $Y \times \mathbb P^{m-n}$ and $X$ with $\deg (Z \to Y \times \mathbb P^{m-n})\le d-2$. So $\deg (Z \to Y \times \mathbb P^{m-n})<\stabirr(X)$, provided $X$ is very general and $d \ge 2m+2$.
\end{example}
\begin{remark}
    Using the same technique one can construct an example even when $f$ is non-constant but non necessarily generically finite onto its image. In this way one can find an example of a very general hypersurface $X\subset \mathbb P^{n+1}$ of degree $d\ge 2m+2$ and a smooth hypersurface $Y\subset \mathbb P^{n+1}$ of degree $e\ge 2n+2$ such that there is a correspondence $Z\subset X \times Y$ with 
    \[
    \deg(Z \to X) = d-2.
    \]
    This shows that the ``very general" hypothesis is needed Theorem A of \cite{lazarsfeld2021measures}.
\end{remark}
\begin{example}\label{bad curve and ab surface}
Let $C$ be a smooth genus $2$ curve and $JC$ its Jacobian. Let also $\varphi: C \to \mathbb P^1$ be the hyperelliptic morphism and take $Z:=C\times C$. We have maps $\id \times \varphi : Z \to C \times \mathbb P^1$ and $Z \to \Sym^2C \sim_\text{birat} JC$. Both these maps have degree $2$, making $Z$ a correspondence between $C \times \mathbb P^1$ and $JC$ with $\deg Z=4$. on the other hand $\gon (C)=2$ and $\irr (JC) \ge \stabirr(JC)\ge 3$ for very general $[C]\in \mathcal M_2$ (see lemma \ref{bound on stab.irr(A)}). So Theorem \ref{curves vs abelian surface} fails if we drop the ``very general" hypothesis.
\end{example}

For the next lemma we recall a definition introduced by Yang \cite{yang2019irrationality}. For any projective variety $X$ of dimension $n$ its \textit{correspondence degree} is defined to be to be
    \[
    \corr (X):= \min \{ \deg (Z \to \mathbb P^n ) \ | \ \text{$Z$ is a correspondence between $X$ and $\mathbb P^n$}\}.
    \]
\begin{lemma}\label{bound on stab.irr(A)} Let $A$ be any abelian surface. Then 
\[
\corr (A)=3 \text{ and } \stabirr (A) \ge 3.
\]
\end{lemma}
\begin{proof}
    The second statement follows from the first since $\stabirr \ge \corr$. To prove $\corr(A)=3$, a similiar argument as the one used to prove $\irr(A) \ge 3$ works to show $\corr (A) \ge 3$ (see \cite{chen}, Lemma 3.5). In brief, let $Z\subset \mathbb P^2 \times A$ be a correspondence with projections $a:Z \to \mathbb P^2, b:Z \to A$ and assume $\deg a=2$. Since all maps $\mathbb P^2 \dashrightarrow A$ are constant, so it is the map
    \[
        \mathbb P^2  \dashrightarrow  A,
        t \mapsto \sum_{u\in a^{-1}(t)}u.
    \]
    Up to translating $Z$, we can assume $\sum_{u\in a^{-1}(t)}u=0$, for general $t\in \mathbb P^2$. Then $\mathbb P^2$ must dominate the Kümmer surface $K(A)$ of $A$, which gives a contradiction, since $K(A)$ is not generically covered by rational curves. 
    To show $\corr (A) \le 3$, notice that, being $A$ very general, it can be dominated by a $(1,2)$-polarized abelian surface $(A',L')$ with $L'$ represented by a smooth curve $C\in |L'|$. Then, by Theorem 2 of \cite{yoshihara1996degree}, $\irr (A')=3$ so $\corr(A) \le 3$.
\end{proof}

The next example shows that for every surface $S$ and for a very general abelian surface $A$, there is always a correspondence $Z\subset S \times A$ such that 
\[
\deg (Z \to S)=3.
\]

\begin{example}\label{ strong princ og gen fails for C, A}
    Let $S$ be any surface and $A$ a very general abelian surface. Recall that for every $r>0$, the \textit{generalized Kümmer variety} $K_r(A)$ is defined to be the fiber over $0_A$ of the summation map $\text{Hilb}^{r+1}(A) \to A$. 
    
    Now fix a rational surface $\Sigma \subset K_2(A)$. Such a $\Sigma$ can be constructed using the same technique as Lemma \ref{bound on stab.irr(A)}: let $\iota:A' \to A$ be an isogeny, where $A'$ is a $(1,2)$-polarized abelian surface containing a smooth curve in the complete linear system of the polarization and let $\psi:A' \dashrightarrow \mathbb P^2$ be a degree $3$ dominant map, which it exists by \cite{yoshihara1996degree}; then take $\Sigma$ as the image 
    \[
    \Sigma:=\overline{\im \Bigg\{ 
        \mathbb P^2 \dashrightarrow \text{Hilb}^3 (A), \ t \mapsto \iota ( \psi^{-1}(t))
    \Bigg\}}\subseteq \text{Hilb}^3(A).
    \]
    $\Sigma$ is rational by construction, being dominated by $\mathbb P^2$. Up to translating $\psi$ we have $\Sigma \subset K_2(A)$. Consider the incidence correspondence
    \[
    I:=\{(\xi,z)\in \text{Hilb}^3(A) \times A \ | \ \text{$\xi$ is supported at $z$}\}.
    \]
    $I$ is irreducible since it is the closure of the image in $\text{Hilb}^3(A) \times A$ of the diagonal 
    $\Delta_{14}\subset A^3 \times A$.
    Take now a dominant map $\varphi:S \dashrightarrow \Sigma$ and consider the correspondence
    \[
    \begin{aligned}
    Z&:=\overline{(\varphi \times \text{id}_A)^{-1}(I)}=\\
    &=\overline{\{(x,u) \in S \times A \ | \ \varphi(x,t) \text{ is supported at }u\}}.
    \end{aligned}
    \]
    $Z$ is irreducible. Moreover, by construction 
    \begin{equation}\label{small deg corr}
    \deg (Z \to S) =3.
    \end{equation}
    
\end{example}

Since the very general $(1,d)$-polarized abelian surface with $d\nmid 6$ has degree of irrationality $4$ (see \cite{martin}), we deduce:
\begin{proposition}\label{p strong principle of generality fails for C,A}
    Let $S$ be any surface and $A$ a very general $(1,d)$-polarized abelian surface with $d\nmid 6$. Then, there is a correspondence $Z\subset S \times A$ between $S$ and $A$ such that
    \[
    \deg (Z \to S)<\irr (A).
    \]
\end{proposition}

\begin{remark}
    Let $A,A', \iota, \varphi$ as in Example \ref{ strong princ og gen fails for C, A}. Then similarly as what we have already done, one can look at the rational map
    \[
    F_\varphi: \mathbb P^2 \dashrightarrow K_2(A).
    \]
    This defines a rational surface $\Sigma$ in $K_2(A)$ by
    \[
    \Sigma:= \overline{F_\varphi(\mathbb P^2)}.
    \]
    Clearly $\Sigma$ is an example of a $C^2$-subvariety of $K_2(A)$ (see Definition \ref{C^i-subvariety}), being a constant cycle subvariety. Recall
    \begin{definition}\label{C^i-subvariety} Let $X$ be a projective variety of dimension $2n$.
    \begin{itemize}
        \item Given $x\in X$, we denote $O_x \subset X$ the set of points of $X$ which are rationally equivalent to $x$;
        \item (\cite{Voisin remarks and questions}) For $0 \le i \le 2n$, a $C^i$\textit{-subvariety} is a subvariety $Y\subset X$ of X of dimension $2n-i$ such that for every $y\in Y$ one has $\dim O_y \ge i$.\footnote{Since $O_y$ is a countable union of Zariski-closed subsets of $X$ with $\dim O_y$ we mean the maximal dimension of a Zariski-closed subset contained in $O_y$.}
    \end{itemize}
    \end{definition}
    More interesting examples of $C^i$-varieties in generalized Kümmers have been found by Lin in \cite{Lin}.
\end{remark}

\section{The stable correspondence degree of two hypersurfaces}\label{section 3}
We devote this section to the proof of Theorem \ref{curves & hypersf}. We will adopt the following notation.
\subsubsection*{Notation}: Let $Z$ and $V$ be projective varieties with a dominant map $\varphi: Z \to V$. Let also $W\subset V$ be a subvariety. We denote
\[
Z|_W:=\varphi^{-1}(W).
\]
\begin{proof}[Proof of Theorem \ref{curves & hypersf}]
Set $d:= \deg X$ and let 
\[
\xymatrix{
        & Z \ar[dl]_a \ar[dr]^b& \\
        X\times \mathbb P^r & & Y \times \mathbb P^s
}
\]
be a correspondence between $X \times \mathbb P^r$ and $Y \times \mathbb P^s$ and assume $\deg b \le \text{stab.irr}(X)$, aiming to show $\deg b = \stabirr (X)$. Hence, $\deg b\le d-1$ by Theorem B of \cite{yang2019irrationality}. First we want to reduce to the case $r=0$. To this end, fix $n_1,s_1 \ge 0$ with $n_1+s_1=m$ (recall that $m=\dim X$ and $n=\dim Y$ by hypothesis) such that for $Y_1\subset Y$ a general $n_1$-dimensional linear section of $Y$ and $P_1\subset \mathbb P^s$ a general $s_1$-dimensional linear subspace of $\mathbb P^s$, $Z|_{Y_1 \times P_1}$ dominates $X$. Then, 
\[
Z_1:=\pr_{X\times Y_1 \times P_1}(Z|_{Y_1 \times P_1})
\]
with projections $a_1:Z_1\to X$ and $b_1:Z_1 \to Y_1\times P_1$, is a correspondence between $X$ and $Y_1 \times P_1$. We observe that 
\[
b_1^{-1}(y,v)=\pr_X(b^{-1}(y,v)). 
\]
Therefore $\deg b_1 \le d-2$, so from now on we assume $r=0$. 

By Example \ref{ex hypersf vs hypersf}, for a general $(y,v) \in Y \times \mathbb P^{n+1}$ the $0$-cycle \[
X_{(y,v)}:=b^{-1}(y,v) \subset X
\]
is not separated by $K_X$. Moreover $\deg X_{(y,v)} \le d-1$, hence by example \ref{ex hypersf vs hypersf}, $X_{(y,v)}$ is made up of at least $d-m$ collinear points (see also \cite{bastianelli2011gonality}). Moreover, $X_{(y,v)}$ spans a line, call it $\ell_{(y,v)}\subset \mathbb P^{m+1}$. Set $\delta:= \deg X_{(y,v)}= \deg  b^{-1}(y,v)$, we can write
\[
(X\cdot \ell_{(y,v)})=X_{(y,v)}+F_{(y,v)}
\]
for some $0$-cycle $F_{(y,v)}$ of degree $d-\delta$.

We are now going to borrow ideas from \cite{BDELU} and \cite{lazarsfeld2021measures}, aiming to show $F_{(y,v)}$ is a point.  

Turning to details, set $B:=Y \times \mathbb P^s$ and let $\mathbb G$ be the Grassmanian of lines in $\Prj^{m+1}$. 
Resolving the singularities of the map 
\[
B \dashrightarrow \mathbb G, (y,v) \mapsto \ell_{(y,v)}
\]
we get a diagram
\begin{equation}
\xymatrix{\mathbb P^{m+1} & W' \ar@{->>}[l]_\mu \ar[d]_\pi \ar[r] & W \ar[d]^{\Prj^1-\text{bundle}} \\
     & B' \ar[r] & \mathbb G}
\end{equation}
where $W\to \mathbb G$ is the tautological projective bundle over $\mathbb G$ and $W'$ its pull-back to $B'$. Note that $\mu$ is generically finite.
Set now $X^*:= \mu^*(X)$, so that
\[
X^*=X'+F
\]
where $X'$ is a reduced and irreducible divisor of relative degree $\delta$ over $B'$. Note that $X'$ represents birationally the correspondence $Z$. If we prove $\deg (F \to B') \le 1$ we will be done. 

\subsubsection*{Claim} Let $V\subset F$ be an irreducible component of $F$ dominating $B'$. Then $\mu(V)$ is a point. 

First we show that the claim implies $\deg (F \to B')\le 1$. Suppose $V_1,\dots, V_N \subset F$ are the irreducible components dominating $B'$. By the claim, $\mu(V_i)$ is a single point for all $i$, call it $p_i \in X$. 
Since $\ell_{(y,v)}$ is a line through $p_1, \dots, p_N$, for $(y,v)\in Y \times \mathbb P^s$ general, 
and  the $0$-cycles 
\[
\{(\ell_{(y,v)} \cdot X) \ | \ (y,v)\in Y \times \mathbb P^s\}
\]
span $X$, it follows that we must have $p_1=\dots =p_N$. So in particular $\mu(F)_\text{red}=p_1$. Therefore, if $\deg (F \to B') \ge 2$, then $\ell_{(y,v)}$ must be a line tangent to $X$ in $p_1$. But then, if $H$ is the tangent hyperplane to $X$ at $p_1$, we have
\[
X=\mu(X')\subseteq \overline{ \bigcup_{(y,v) \in Y \times \mathbb P^s}\ell_{(y,v)}}\subseteq H.
\]
This is a contradiction, so we must have $\deg (F \to B')=1$ and Theorem \ref{curves & hypersf} follows.

As for the claim, let $V\subset F$ be an irreducible component dominating $B'$. Set 
\[
T:= \mu (V)\subset X.
\]
We claim first that $T \subsetneq X$; in fact if $T=X$ then we could construct a new correspondence between $X$ and $Y \times \mathbb P^s$ having degree over $Y \times \mathbb P^s$ equal to $\deg (V \to B')\le m<d-m$, which we showed before cannot happen. 
Now let $t:=\dim T$. 
We want to prove $t=0$. View $V$ as reduced irreducible variety of dimension $n$. After desingularizing we get the following morphisms of smooth varieties
\[
\xymatrix{X & T' \ar[l] & V' \ar[l]_{\mu'} \ar[d]_{\pi'} \\
& & B'}
\]
Denote by $\sigma \ge 0$ be the maximal integer with the following property:
if $\nu:=t-\sigma$ then,
for 
\[
Y_L \subset Y, P_L\subset \mathbb P^s
\]
   respectively $\nu$-dimensional general linear subspace of $Y$ and $\sigma$-dimensional general linear subspace of $\mathbb P^s$, we have that
\[
V_L:=V'|_{Y_L \times P_L}
\]
dominates $T'$.
Setting then $B_L:= Y_L \times P_L$ and desingularizing we get the following commutative diagram of smooth varieties:
 \begin{equation}\label{cursed diagram}
\xymatrix{X & T' \ar[l] & V' \ar[l]_{\mu'} \ar[d]_{\pi'_{V'}} & V_L' \ar[d] \ar[l] \\
& & B' & B_L' \ar[l]^j  }
\end{equation}
with $V_L' \to V'$ and $B_L'\to B'$ birational onto their images.
Set $e:= \deg (V' \to B')$. We distinguish now two cases. 

\subsubsection*{Case 1} Assume $\sigma>0$. Then ${V_L'}_*:H^{t,0}(T') \to H^{t,0}(B_L')$ vanishes, since $B_L'$ is ruled; therefore $K_{T'}$ does not separate the fibers of $V' \to B'$, hence $K_{T'}$ does not satisfy $\BVA{e-1}$.
Now, if $t>0$, by result of Ein and Voisin \cite{ein1988subvarieties} \cite{voisin1996conjecture}, $K_{T'}$ satisfies  $(\text{BVA})_{d+t-2m-2}$ (see \cite{BDELU}, Proposition 3.8), so one must have 
\[
e \ge d+t-2m.
\]
On the other hand, 
computations with the canonical bundles give $e(m-t)\le m$.\footnote{As in \cite{BDELU} proof of Theorem C, page 2383: by standard computations $(K_{W'/\mathbb P^{m+1}}\cdot F)=m$, for general $F$ fiber of $W' \to B$; this together with $\text{ord}_V(K_{W/\mathbb P^{m+1}}) \ge m-t$ (see \cite{BDELU} Appendix A, corollary A.6 for a proof) gives the desired inequality.} So, as in \cite{BDELU} and \cite{lazarsfeld2021measures} we have 
\[
d-2m \le \frac m{m-t}-t.
\]
Now the expression
\[
\frac m{m-t}-t 
\]
is always $\le 1$ if $1 \le t \le m-1$, which forces $d\le 2m+1$. Therefore $t=0$ and we are done.

\subsubsection*{Case 2} Suppose that $ \sigma=0$, so in particular $\nu =t$.  Then the $0$-cycle
\[
V'_{(y,v)}:=\mu'_*{\pi'}_{V'}^*(y,v) 
\]
depends only on the $Y$-component,\footnote{Here with abuse of notation we identify  $\pi'_{V'}$ with $V' \to B' \to B$. } i.e.: for a general $y\in Y$ and for general $v,w\in \mathbb P^s$ one has
    \[
    V'_{(y,v)}=V'_{(y,w)}.
    \]
Suppose $\#V'_{(y,v)} \ge 2$, for $(y,v)\in Y \times \mathbb P^s$ general.
Then $\ell_{(y,v)}$ deepends only on the $Y$-component, being $\ell_{(y,v)}=\langle V'_{(y,v)}\rangle$. Therefore also $X_{(y,v)}$ depends only on the $Y$-component. Hence $s=0$, because $X$ is covered by $X_{(y,v)}$  as $y$ moves in $Y$, for some general fixed $v\in \mathbb P^s$. Then the theorem follows from Theorem A of \cite{lazarsfeld2021measures}. In fact Lazarsfeld and Martin (see \cite{lazarsfeld2021measures} page 46) prove that if $s=0$ and $\deg b \le d-1$, then $\deg b = d-1$ and $T \subset \mu(F)= \pt$.

Assume now $\# V'_{(y,v)}=1$, for $(y,v)\in Y \times \mathbb P^s$ general. Then, fixed a general $v_0\in \mathbb P^s$, we can define the rational map
\[
f:Y \dashrightarrow X, y \mapsto  V'_{(y,v_0)}
\]
By Theorem 4.2 of \cite{families}, there are no non-constant rational maps $Y \dashrightarrow X$. So $f$ needs to be constant and, since $T=\overline{f(Y)}$, we are done.
\end{proof}
\begin{remark}
    Let $X\subset \mathbb P^{m+1}$ be a very general hypersurface of degree $\ge 2m+2$ and $Y$ a very general curve of genus $g$. Take $Z$ as in the Set-Up. Then by Example \ref{ex curve vs hypersf}, $\deg a \ge \gon (Y)$. 
    Moreover, the same argument used in the proof just completed together with Theorem A of \cite{families} shows $\deg b \ge \deg X-1$. So we find
\end{remark}
\begin{theorem}\label{curves vs hypersf 2}
    Let $X \subset \mathbb P^{m+1}$ and $C$ be respectively a very general hypersurface of degree $\ge 2m+2$ and a very general curve of genus $g$. Then
    \[
    \stabirr(X,C)= \stabirr(X) \cdot \gon (C).
    \]
\end{theorem}

\section{Curves and Abelian surfaces}\label{section 4}

In this section we prove Theorem \ref{curves vs abelian surface}. We also state and prove a partial result regarding the stable correspondence degree of a very general curve and an abelian surface.

We start with some terminology.

\begin{definition}
    Let $X$ be an $n$-dimensional smooth projective variety. Consider an integer $\delta >1$. Given a $n$-dimensional projective variety $W\subset \text{Hilb}^\delta(X)$, we define its \textit{multiplicity} as 
    \[
    \mult (W):=  W \cdot (x+ \text{Hilb}^{\delta-1}(X)) \in \mathbb Z,
    \]
    for a general $x\in X.$
\end{definition}
We have the following simple proposition.
\begin{proposition}\label{irr con mult}
    Let $X$ as in the previous definition and $c >0$. Then
    \[
    \irr (X) \le c
    \iff 
    \begin{aligned}
        &\exists \delta \le c \text{ and a rational variety }W\subset \textnormal{Hilb}^\delta(X) \text{ with }\\
        &\dim W=n \text{ and } \mult (W)=1.
    \end{aligned}
    \]
\end{proposition}
\begin{proof} If $\irr (X) \le c$, then there is  a dominant rational map $\varphi: X \dashrightarrow \mathbb P^n$ with $\delta := \deg \varphi \le c$. This defines an injective rational map
\[
F_\varphi: \mathbb P^n \dashrightarrow \text{Hilb}^\delta(X), t \mapsto x_1+ \dots +x_\delta 
\]
where $\varphi^{-1}(t)=\{x_1, \dots, x_\delta\}$. Then $W:= \overline{\im F_\varphi}$ is an $n$-dimensional rational variety. Let now $x\in X$ general and let 
\[
x_1+\dots +x_\delta, y_1+ \dots +y_\delta \in W\cdot (x+\Hilb^{\delta-1}(X)).
\]
Then
\[
\{x_1, \dots ,x_\delta\}=\varphi^{-1}(\varphi(z))=\{y_1, \dots, y_\delta\} \implies x_1+\dots +x_\delta = y_1+\dots +y_\delta.
\]
This shows $\mult (W)=1$.

Let now $W\subset \Hilb^\delta (X)$ be a rational $n$-dimensional variety with $\mult (W)=1$, for some $\delta \le c$. 
Consider then the map
\[
\begin{aligned}
\psi=\psi_W :&X \dashrightarrow W \sim_\text{birat.} \mathbb P^n,\\
&x \mapsto W\cdot (x+\Hilb^{\delta-1}(X)).
\end{aligned}
\]
Let $\xi =w_1+ \dots +w_\delta \in W$ be a general point. Then $\psi(w_1)=\xi$, so $\psi$ is dominant. On the other hand
\[
\psi^{-1}(\xi)\subseteq \{w_1, \dots, w_\delta\}
\]
so $\deg \psi \le \delta \le c$.
\end{proof}

We are now ready to prove Theorem \ref{curves vs abelian surface}.

\begin{proof}[Proof of Theorem \ref{curves vs abelian surface}]
Let $Z\subset (C\times \mathbb P^1) \times A$ be a correspondence, with maps
\[
    \xymatrix{
        & Z \ar[dl]_a \ar[dr]^b& \\
        C\times \mathbb P^1 & & A
}
\]
 We have already proven that $\deg b \ge \text{gon}(C)$, see Example \ref{ex curve ab surf}. Therefore, since $\irr (A) =3 \text{ or }4$ (see \cite{chen} and \cite{martin}), it is enough to show that $\deg a >2$ and that if $\deg Z<4 \cdot \text{gon} (C)$ then $\irr (A)=3$. 
 
 Since for very general $C$, all maps $C \times \mathbb P^1 \dashrightarrow A$ are constant, so is the map
\[
C \times \mathbb P^1 \dashrightarrow A, (x,t) \mapsto \sum_{u \in a^{-1}(x,t)}u.
\]
Therefore, up to translating $Z$ we can assume such a map to be constantly equal to $0=0_A$. 

Assume $\deg a=2$. Then we have a dominant map $C \times \mathbb P^1 \dashrightarrow K(A)$, where $K(A)$ is the Kümmer surface of $A$.
This is a contradiction since $h^0(K(A),\omega_{K(A)})>0$.

Assume now $\deg Z< 4\cdot \gon (C)$, hence $\deg a=3$ and $\deg b < \frac 43 \gon (C)$. Then as before we get a map
\[
\varphi=\varphi_Z: C \times \mathbb P^1 \dashrightarrow K_2(A)
\]
where $K_r(A)$ is the generalized Kümmer. Such a map must be generically finite onto its image. Therefore \[
\Sigma:=\overline{\varphi (C \times \mathbb P^1)}
\]
is a 
(possibly singular) surface in $K_2(A)$. By Proposition \ref{ruled surfaces in 4-folds w/ 2-forms} or Example \ref{ruled surfaces in HK 4-folds} (see below) it follows that $\Sigma$ is a rational surface.
So if we show that $\mult (\Sigma)=1$, by Proposition \ref{irr con mult}, we will have $\irr (A)\le 3$ and we will be done.

As for the proof of $\mult (\Sigma)=1$, notice first that for general $u+v+w\in \im \varphi$ we have
\begin{equation}\label{fibers of phi}
\varphi^{-1}(u+v+w)\subset b^{-1}(u).
\end{equation}
In fact, if $(x,t)$ belongs to the left hand side of (\ref{fibers of phi}), then $(x,t,u)\in Z$, hence $(x,t)\in b^{-1}(u)$. Consider now a general $u\in A$ and let
\[
u+v+w,u+v'+w' \in \Sigma \cdot (u + \Hilb^2(A)).
\]
Then, by (\ref{fibers of phi}) we have
\[
\varphi^{-1}(u+v+w), \varphi^{-1}(u+v'+w')\subset b^{-1}(u).
\]
Since $\deg \varphi \ge \gon (C)$ and $\# b^{-1}(u)<\frac 43 \gon (C)$, we have that $\varphi^{-1}(u+v+w) \text{ and } \varphi^{-1}(u+v'+w')$ 
must intersect, hence \[
u+v+w=u+v'+w'\text{ as a points of }\Hilb^3(A).
\]
So $\mult (\Sigma)=1$ and we are done. 
\end{proof}

\begin{remark}
Let $Z\subset (C\times \mathbb P^{r+1}) \times (A \times \mathbb P^r)$ be a correspondence  between $C \times \mathbb P^{r+1}$ and $A\times \mathbb P^r$ for some $r\ge 0$, with projections $a:Z \to C \times \mathbb P^{r+1}$ and $ b :Z \to A \times \mathbb P^r$. The same technique used to prove Theorem \ref{curves vs abelian surface} works to show $\deg b \ge \gon (C)$ and $\deg a \ge 3$. Moreover if $A$ is $(1,d)$-polarized for $d=2$ or $6$, $\irr (A)=3$, see \cite{yoshihara1996degree} and \cite{moretti2023polarized}. Therefore by equation (\ref{obvioous bound}) in the introduction we have $\stabcorr(C,A)=3 \cdot \gon (C)$. Summing everything together we have
\end{remark}
\begin{proposition}\label{stab.corr(C,A)} Let $A$ a very general abelian surface. Then for a very general curve $C$ we have
\[
\stabcorr (C,A) \ge 3 \cdot \gon (C)
\]
and equality holds if $A$ is $(1,d)$-polarized with $d=2$ or $6$.
\end{proposition}

We close the section by proving results that have appeared in the proof of Theorem \ref{curves vs abelian surface}. 

We start with Proposition \ref{ruled surfaces in 4-folds w/ 2-forms}. We would like to thank Olivier Martin for suggesting the following approach, which lead to simpler argument than the one originally found by the author.

\begin{proposition}\label{ruled surfaces in 4-folds w/ 2-forms}
    Let $X$ be a smooth irreducible, projective $4$-dimensional variety, having a holomorphic $2$-form of maximal rank that does not drop rank on a $3$-fold i.e.: there is $\omega \in H^0(X, \Omega^2_X)$ such that $\dim (\{x\in X \ | \ \textnormal{rank}(\omega_x)<4\}) <3$.\footnote{By "$\rk(\omega_x)$" we mean the rank of $\omega_x$ as a skew-symmetric linear $2$-form on $T_xX$, for every $x\in X$.}  
    Let $C$ be a very general curve of genus $g> 1$. Then any rational map
    \[
    C \times \mathbb P^1 \dashrightarrow X
    \]
    which is generically finite onto its image, factors through a rational surface.
\end{proposition}
\begin{example}\label{ruled surfaces in HK 4-folds}
    In particular from Proposition \ref{ruled surfaces in 4-folds w/ 2-forms} we have that if $X$ is a HyperKähler $4$-fold, for a very general curve $C$ of genus $g>1$, all rational maps $C \times \mathbb P^1 \dashrightarrow X$ which are generically finite onto their images factor through $\mathbb P^2$.
\end{example}
\begin{example}
    Note that there are complex $4$-folds satisfying the hypothesis of Proposition \ref{ruled surfaces in 4-folds w/ 2-forms} but containing a birational copy of every genus $g$ curve, for all $g \ge 0$. 
    For example, if $A$ is an abelian surface, $X:=K_2(A)$ contains a rational surface (see for example Example \ref{ strong princ og gen fails for C, A}). So in particular, given a curve $C$ there is always a morphism $C \to X$ which is birational onto its image.
\end{example}
\begin{proof}[Proof of Proposition \ref{ruled surfaces in 4-folds w/ 2-forms}]
Let $\Sigma:=\overline{f(C \times \mathbb P^1)}$. Then $\Sigma$ is a surface with $\text{kod}(\Sigma)= -\infty$. Therefore $\Sigma$ is birationally equivalent to to $D \times \mathbb P^1$ for some smooth curve $D$. If $g(D)=0$ there is nothing to prove, so assume $g(D)>0$. 

There is a dominant map $C \to D$. Since a very general curve $C$ of genus $>0$ can only dominate $\mathbb P^1$ or itself, we have $D\cong C$. This shows that for a very general $[C] \in \mathcal M_g$, there is a rational map
\[
C \times \mathbb P^1 \dashrightarrow X
\]
which is birational onto its image. We want to show this cannot happen. Fix a copy $\mathbb P$ of $
\mathbb P^1$, say $\mathbb P := \mathbb P (\mathbb C^2)$. 
By Proposition \ref{families of rational maps}, we can construct a family 
\[
\mathcal C \to B
\]
of genus $g$ curves over an irreducible and reduced quasi-projective variety $B$, whose natural map $B \to \mathcal M_g$ is dominant and generically finite, together with a rational map
\[
F: \mathcal C \times \mathbb P \dashrightarrow X
\]
such that, for every $b \in B$, the restriction
\[
F_b: C_b \times \mathbb P \dashrightarrow X
\]
is birational onto its image. Let $V \subseteq X$ be the closure of the image of $F$. Then $\dim V=3$. In fact, being $F_b$ birational for every $b\in B$, it must be $\dim V>2$. On the other hand, $V\subsetneq X$ since $X$ is not covered by rational curves, being $h^0(\Omega_X^2)>0$.  

By \cite{rational curves} Theorem 5.2, the maximal rationally connected fibration exists for $V$, i.e.: there is an Zariski-open subset $V^0 \subseteq V$, an irreducible quasi-projective variety $Z^0$ and a dominant morphism
\[
\pi:V^0 \to Z^0
\]
whose fibers are rationally connected\footnote{This is because in characteristic zero a smooth and rationally chain connected variety is rationally connected by \cite{rational curves} Theorem 3.10.3}. 

By a theorem of Graber, Harris and Starr (\cite{Graber-Harris-Starr} Corollary 1.4), the base of an maximal rationally connected-fibration cannot be uniruled. So $Z^0$ is not uniruled.

Since $\dim V=3$, we have that $\omega$ has maximal rank at the general point of $V$. Therefore $\omega$ defines a holomorphic $2$-form on $V$, which is non-zero by linear algebra. So $V$ is not rationally connected, hence $\dim Z^0>0$. 

Assume $\dim Z^0=1$. Let $Z$ be the normalization of a projective closure of $Z^0$. Then $Z$ is covered by all genus $g$ curves, so it must be $Z \cong \mathbb P^1$. So we reached a contradiction with the fact that $Z$ cannot be uniruled.

So $\dim Z^0=2$. Let $Z$ be a projective closure of $Z^0$ and 
\[
\mu: Z \dashrightarrow S
\]
be a birational map to a smooth minimal surface $S$. Let $\Phi:= \mu \circ  \pi \circ F$. Since, for a very general $z\in Z$, all rational curves in $V$ intersecting $\pi^{-1}(z)$ are contained in $\pi^{-1}(z)$ (see \cite{rational curves} Complement 5.2.1 for a proof), we have that $\Phi|_{C_b \times \mathbb P}$ must be constant in the $\mathbb P$-entry. This defines a morphism
\[
\Phi_b: C_b \to Z,
\]
for a very general $b\in B$.
The images $\Phi_b(C_b)$ cover $S$. Since a very general curve of genus $g$ can only dominate curves birationally equivalent to itself or to $\mathbb P^1$, we have two cases.

Case 1: $\Phi_b(C_b)$ is rational, for a very general $b\in B$. This means that $S$ is uniruled, which is a contradiction. 

Case 2: $\Phi_b(C_b)$ is birational to $C_b$ for a very general $b\in B$. Then $S$ needs to be rational, by Lemma \ref{final lemma}. This is again a contradiction and we are done.

\end{proof}

\begin{lemma}\label{final lemma}
    Fix an integer $g> 1$. Let $S$ be a smooth minimal surface. Suppose that for a very general curve $C$ of genus $g$ there is a map $C \to S$ which is birational onto its image. Then $S$ is rational.
\end{lemma}
\begin{proof} It is enough to show that $S$ is ruled. In fact then $S$ is birational to $X \times \mathbb P^1$, for some smooth curve $X$. Then, $X$ is dominated by all genus $g$ curves, so it must be $X\cong \mathbb P^1$. 

 Suppose $S$ is not ruled. Then $\text{kod}(S) \ge 0$. As before, using Proposition \ref{families of rational maps}, we can construct a family of genus $g$ curves $ \pi: \mathcal C \to B$ over an irreducible quasi-projective variety, whose natural map $B \to \mathcal M_g$ is dominant and generically finite, together with a rational map
 \[
 F: \mathcal C\dashrightarrow S
 \]
 with $F_b:C_b \to S$ birational onto its image, for a very general $b\in B$. The indeterminacy locus of $F$ has codimension $\ge 2$, so it does not intersect all the fibers of $\pi$. Therefore, up to replace $B$ with a Zariski-open subset, we can assume $F$ is a morphism. This shows that, for $b\in B$ general, the space of deformations with fixed target of the morphism $F_b:C_b \to S$ has dimension $ \ge 3g-3$.

 Now, a result by Lee and Pirola (Proposition 1.2 of \cite{subfield}) gives that if $U$ is the Kuranishi space of deformations of a birational morphism from a genus $g$ curve to a surface with nef canonical bundle, then
 \[
 \dim U \le g.
 \]
 Now, since 
 $S$ is minimal and $\text{kod}(S) \ge 0$, we have that $K_S$ is nef (see Example 1.4.18 of \cite{positivity in algebraic geometry}). 
 So we have 
 \[
 g \ge 3g-3,
 \]
 which contradicts $g \in \mathbb Z_{>1}$. So $S$ must be ruled.
\end{proof}

Finally we prove the following statement used in the proofs of Proposition \ref{ruled surfaces in 4-folds w/ 2-forms} and Lemma \ref{final lemma}.

\begin{proposition}\label{families of rational maps} 
Let $V$ and $\Gamma$ be a projective varieties. Assume that for a very general curve $C$ of genus $g\ge 0$ there is a rational map
\[
C \times \Gamma \dashrightarrow V
\]
which is generically finite of degree $d$ onto its image. Then there is a family of genus $g$ curves
\[
\mathcal Y \to B
\]
over an irreducible quasi projective variety $B$ such that the natural map $B \to \mathcal M_g$ is dominant, together with a rational map
\[
F: \mathcal Y \times \Gamma \dashrightarrow V
\]
such that for a general $b\in B$ the restriction
\[
F_b: X_b \times \Gamma \dashrightarrow V
\]
is generically finite of degree $d$ onto its image.
\end{proposition}
\begin{proof}
    Let $e: M \to \mathcal M_g$ be a quasi projective generically finite cover of the moduli space of curves, with a family of curves $\mathcal C \to M$ such that for a general $m \in M$ one has that $C_m$ is a smooth genus $g$ curve which belongs to the isomorphism class $e(m)\in \mathcal M_g$ (take $\mathcal M_g$ to be the coarse moduli space in the case $g=1)$. \footnote{This could be achieved in many ways. One way could be seeing $\mathcal M_g$ as a GIT quotient $H \to \mathcal M_g$ of some Hilbert scheme $H$ in $\mathbb P^N$, for some $N=N_g\gg g$. $H$ admits a tautological family $\mathcal H \to H$. Then take $M$ as a suitable linear section of $H$ (recall that $H$ is quasi projective) and $\mathcal C:= \mathcal H_M$.}
    
By hypothesis, for general $m \in M$ there is a map $C_m \times \Gamma \dashrightarrow X$ which is generically finite of degree $d$ onto its image. 

We borrow ideas from the proof of Theorem 1 of \cite{Voisin Nori}. 
By properties of Chow varieties or relative Hilbert schemes there are countably many varieties $\mathcal B_k$ parametrizing pairs $(C_m,\varphi)$ where $m \in M$ and $\varphi:C_m \times \Gamma \dashrightarrow V$ is generically finite of degree $d$ onto its image. In fact such maps can be identified with their graph, which is some subvariety of dimension $\dim (\Gamma)+1$ in $C_m\times \Gamma \times X$, moving in some parameter space.

Let now 
\[
\mathcal Y_k:= \mathcal B_k \times_M \mathcal C.
\]
There is a natural rational map 
\[
F_k: \mathcal Y_k \times \Gamma \dashrightarrow X.
\]
By Baire Cathegory Theorem $\mathcal B_{k_0}$ must dominate $M$ for some $k_0$.
After slicing if needed $\mathcal B_{k_0}$ by a suitable (general) linear section and restricting our attention to an irreducible component, we get a family 
\[
\mathcal Y \to  B
\]
of curves of genus $g$ over an irreducible base, which pulls back from $\mathcal C \to M$ via a generically finite dominant morphism $B \to M$ and having a rational map 
\[
F: \mathcal Y \times \Gamma \dashrightarrow V
\]
such that 
\[
F_b:Y_b \times \Gamma \dashrightarrow X 
\]
is generically finite of degree $d$ onto its image, for $b\in B$ general. This is the desired family.
\end{proof}

 %
 %
 %
 %

 \end{document}